\documentclass{amsart}
\usepackage{amsrefs, amsmath,amsthm,amssymb,enumerate,mathrsfs}
\usepackage{xcolor}
\usepackage[colorlinks=true,  linkcolor=red, citecolor=blue]{hyperref}
\newtheorem{theorem}{Theorem}[section]
\newtheorem{lemma}[theorem]{Lemma}

\newtheorem{proposition}[theorem]{Proposition}

\newtheorem{remark}[theorem]{Remark}

\numberwithin{equation}{section}

\begin{document}
\title[]{Locally constrained inverse curvature flow and Alexandrov-Fenchel type inequalities in de Sitter space}

\author{Kuicheng Ma}
\address{School of Mathematics and Statistics,
Shandong University of Technology, Zibo 255000, China}
\email{kchm@sdut.edu.cn}





\begin{abstract}
In this paper, we study the behavior of some locally constrained inverse curvature flow in de Sitter space, with initial value any closed spacelike $k$-convex hypersurface satisfying some pinching condition. Assume further the Heintze-Karcher inequality for any closed spacelike mean convex hypersurface in de Sitter space, we derive a class of Alexandrov-Fenchel inequalities.
\end{abstract}

\maketitle


\par

\section{Introduction}
Let $(\mathbb{S}^n, \sigma)$ be the standard spherical space of dimension $n$, here and in the sequel, $n$ is an integer no less than $2$. Consider the product manifold
\begin{align}\label{pm}
(a, b)\times \mathbb{S}^n,
\end{align}
together with its semi-Riemannian structure
\begin{align}\label{as}
\bar g = \epsilon \, dr^2+\lambda^2(r)\sigma,
\end{align}
where $\epsilon$ is either $1$ or $-1$ and $\lambda:(a,b)\subset \mathbb{R}\rightarrow\mathbb{R}$ is some smooth positive function. When $\epsilon=-1$, $(a, b)=(0, +\infty)$ and $\lambda(r)=\cosh(r)$, \eqref{as} recovers the de Sitter space $\mathbb{S}^n_1$, which is the Lorentzian space form of constant sectional curvature $1$. When $\epsilon=1$, \eqref{as} may recover Riemannian space forms. Furthermore, \eqref{as} admits a conformal Killing vector field $V=\overline \nabla \Lambda$ in the sense
\begin{align}\label{ck}
\overline\nabla V=\epsilon \, \lambda'\, \mathcal{I},
\end{align}
where $\overline\nabla$ is the Levi-Civita connection determined by \eqref{as}, $\Lambda$ is a primitive of $\lambda$ and $\mathcal{I}$ is the identical transformation. Notice that $\Lambda$ can be chosen as $\lambda'$ in hyperbolic space as well as in de Sitter space.

Let $x_0: \mathbb{S}^n\rightarrow \mathbb{S}^n_1$ be a spacelike embedding. In other words, the induced structure on $\Sigma_0:=x_0(\mathbb{S}^n)$ is a Riemannian one. Denote by $\nu$ the future-directed timelike unit normal vector field along $x_0$.  Let $h$ and $\mathcal{W}$ be respectively the second fundamental form and the Weingarten transformation with respect to such a normal vector field. Let $(\kappa_1,\cdots,\kappa_n)$ be the principal curvatures and $E_k$ be the normalized $k$-th elementary symmetric function in principal curvatures, and let $E_0=1$. Define
\begin{align}
u=\bar g(V,\nu)\nonumber
\end{align}
and refer to it as the support function occasionally. For an integer $k\in\{1,\cdots,n\}$, a spacelike hypersurface is said to be $k$-convex if $E_k>0$ holds everywhere on it. As well known, $n$-convex is convex and $1$-convex is mean convex. In general, $k$-convexity implies the principal curvatures $(\kappa_1,\cdots, \kappa_n)$ belong to
\begin{align}
\Gamma_k^+=\{(x_1,\cdots, x_n)\in\mathbb{R}^n~|~ E_l(x_1,\cdots, x_n)>0, ~l=1,\cdots, k\},\nonumber
\end{align}
which contains the positive cone
\begin{align}
\Gamma^+=\{(x_1,\cdots, x_n)~|~ x_1>0,\cdots, x_n>0\}.\nonumber
\end{align}
Obviously, $\Gamma^+=\Gamma_n^+\subset \cdots \subset \Gamma_1^+$. For compact hypersurface in the hyperbolic space $\mathbb{H}^{n+1}$, the horo-convexity and the static convexity are often imposed and static convexity is weaker than horo-convexity. It is worth to pointing out that various kinds of convexity conditions, as a particular kind of pinching conditions, are necessary to establish geometric inequalities.

In this paper, we consider in de Sitter space $\mathbb{S}^n_1$ the locally constrained inverse curvature flow
\begin{align}\label{f1}
\frac{\partial }{\partial t} x=\bigg(\frac{u}{\lambda'}-\frac{E_{k-1}}{E_k}\bigg)\nu:=\mathcal{F}\nu
\end{align}
with initial value $x_0$ some compact spacelike $k$-convex hypersurface, as well as its applications. Here and in the sequel, $k\in\{2,\cdots,n\}$ is assumed. For simplicity, let
\begin{align}
\Theta=\frac{u}{\lambda'}.\nonumber
\end{align}
Inspired by the static convexity  for compact hypersurface in the hyperbolic space $\mathbb{H}^{n+1}$ and the dual relationship between convex spacelike hypersurfaces in $\mathbb{S}^n_1$ and those in $\mathbb{H}^{n+1}$, we assume further on the initial value $x_0$ the pinching condition
\begin{align}\label{pc}
\mathcal{W}\leq  \Theta\, \mathcal{I},
\end{align}
which turns out to be preserved along flow \eqref{f1}. Under the assumptions stated above, we show the long-time existence and convergence to coordinate slice of flow \eqref{f1} first. Curvature estimate along extrinsic geometric flows in Lorentzian setting is much more difficult than that in Riemannian setting, if it is not impossible. However, the curvature estimate along flow \eqref{f1} reduces to the preservation of $k$-convexity under the pinching condition \eqref{pc}.

For a compact spacelike hypersurface $\Sigma$ in semi-Riemannian space forms, its quermassintegrals are defined as
\begin{align}
\mathcal{A}_{-1}(\Sigma)&=(n+1)\int_{\hat\Sigma} ~\text{dvol},~\;~ \mathcal{A}_0(\Sigma)=\int_{\Sigma}~ d\mu,\nonumber\\
\mathcal{A}_l(\Sigma)&=\int_{\Sigma} E_l~d\mu-\epsilon K\frac{k}{n-k+2}\mathcal{A}_{k-2} (\Sigma), ~l=1,\cdots, n,\nonumber
\end{align}
where $K$ is the constant sectional curvature of \eqref{as}, $\hat\Sigma$ is the region enclosed by $\Sigma$ and the coordinate slice $\{0\}\times \mathbb{S}^n$, $\text{dvol}$ and $d\mu$ are respectively the volume element and the induced area element. For a compact spacelike hypersurface in space forms, the optimal inequalities between quermassintegrals are usually referred to as the Alexandrov-Fenchel inequalities or the isoperimetric inequalities. Define the weighted curvature integrals of a compact spacelike hypersurface $\Sigma$ in space forms as
\begin{align}
\mathcal{B}_{-1}(\Sigma)=(n+1)\int_{\hat\Sigma} \lambda' ~\text{dvol}+\lambda^{n+1}(0)\omega_n,~\;~\mathcal{B}_l(\Sigma)=\int_{\Sigma} \lambda' E_l~d\mu, ~\;~ l=0,\cdots, n, \nonumber
\end{align}
where $\omega_n$ is the area of $(\mathbb{S}^n,\sigma)$. In case \eqref{as} is Euclidean, $\lambda(r)=r$, $K=0$, then it follows that weighted curvature integrals coincide with quermassintegrals. Hence weighted curvature integrals can be considered as variants of the quermassintegrals, and the optimal inequalities between weighted curvature integrals are thus referred to as Alexandrov-Fenchel type inequalities occasionally. In particular, the optimal inequality between weighted curvature integrals $\mathcal{B}_{-1}$ and $\mathcal{B}_1$ is usually referred to as the Minkowski-type inequality.

Due to the pinching condition \eqref{pc}, the Minkowski formulae and the Newton inequalities, the weighted curvature integral $\mathcal{B}_k(\Sigma_t)$ is monotone increasing along flow \eqref{f1}. To derive geometric inequalities, we assume further the Heintze-Karcher inequality
\begin{align}\label{hk}
\int_{\Sigma} u ~d\mu\leq \int_{\Sigma}\frac{\lambda'}{E_1} ~ d\mu
\end{align}
for any closed spacelike mean convex hypersurface in de Sitter space $\mathbb{S}^n_1$, which is up to now unavailable to the best knowledge of authors. In fact, the Heintze-Karcher inequality \eqref{hk} implies that $\mathcal{B}_{-1}(\Sigma_t)$ is monotone decreasing along flow \eqref{f1}. See \cite{Bre}, \cite{QX} and \cite{WWZ} for more about Heintze-Karcher inequality.

Now the Alexandrov-Fenchel type inequality follows once the long-time existence and convergence of flow \eqref{f1} are at hand. According to the standard regularity theory in PDE, it suffices to derive along flow \eqref{f1} the uniform a priori estimates up to second order, noticing that the short-time existence of flow \eqref{f1} was claimed by Huisken and Polden in \cite{HP}.

In recent years, a great deal of attentions has been directed to optimal geometric inequalities and geometric flows have made great contributions in this direction. Guan and Li \cite{GL} established the full range of Alexandrov-Fenchel inequalities for convex hypersurfaces in Euclidean space $\mathbb{R}^{n+1}$. Later, Wang and Xia \cite{WX} derived the full range of Alexandrov-Fenchel inequalities for horo-convex hypersurfaces in hyperbolic space $\mathbb{H}^{n+1}$, see also \cite{HLW}. For hypersurfaces in spherical space $\mathbb{S}^{n+1}$, Chen and Sun \cite{CS} obtained the Alexandrov-Fenchel inequalities between $\mathcal{A}_{k}$ and $\mathcal{A}_{k-2}$ for every $k\in\{1,\cdots, n\}$. For convex spacelike hypersurfaces in de Sitter space $\mathbb{S}^n_1$, the Alexandrov-Fenchel inequalities are far away from being complete, see \cite{HL-23, LS, M23, M24,JS}. Using the inverse mean curvature flow, Brendle, Hung and Wang \cite{BHW} derived a Minkowski-type inequality for mean convex hypersurfaces in anti de Sitter-Schwarzschild space. Scheuer and Xia \cite{SX} considered a flow of type \eqref{f1} in Riemannian warped products and they derived for mean convex hypersurface in anti de Sitter-Schwarzschild space as well a Minkowski-type inequality, which is slightly different from that one in \cite{BHW}. Later, main results in \cite{SX} were extended by Li and Ma \cite{LM} to de Sitter space $\mathbb{S}^n_1$. In \cite{HLW}, the optimal inequalities which compare the weighted curvature integrals with the quermassintegrals of horo-convex hypersurfaces in hyperbolic space $\mathbb{H}^{n+1}$ were also established and such kind of optimal inequalities are as well referred to as Aexandrov-Fenchel type inequalities there. Very recently, Hu and Li \cite{HL-22} showed for static convex hypersurface in hyperbolic space $\mathbb{H}^{n+1}$ the Alexandrov-Fenchel type inequality between $\mathcal{B}_l$ and $\mathcal{B}_{-1}$ with $l\in\{0,\cdots, n\}$, and they conjectured the full range of Alexandrov-Fenchel type inequalities under static convexity condition. Besides, the authors in \cite{HL-22} recovered the Alexandrov-Fenchel type inequalities shown in \cite{HLW} after replacing the horo-convexity with static convexity. One main tool employed by authors in \cite{HL-22} is a locally constrained mean curvature flow along which the static convexity is preserved.

Main results in this paper are formulated in the following theorems.
\begin{theorem}
Let $k\in\{2,\cdots,n\}$ and $\Sigma_0$ be a closed spacelike hypersurface which is $k$-convex and satisfies the pinching condition \eqref{pc} in de Sitter space $\mathbb{S}^n_1$, then the locally constrained inverse curvature flow \eqref{f1}, with $\Sigma_0$ as its initial value, exists for all positive time and converges to some coordinate slice.
\end{theorem}

\begin{theorem}\label{MT2}
Let $k\in\{2,\cdots,n\}$ and $\Sigma_0$ be a closed spacelike hypersurface which is $k$-convex and meets the pinching condition \eqref{pc} in de Sitter space $\mathbb{S}^n_1$. Assume further the Heintze-Karcher inequality \eqref{hk} for any closed spacelike mean convex hypersurfaces in de Sitter space $\mathbb{S}^n_1$. Then there exist monotone increasing functions $\varphi_{-1}$ and $\varphi_k$ such that
\begin{align}
\mathcal{B}_k(\Sigma_0)\leq \varphi_k\circ \varphi^{-1}_{-1}(\mathcal{B}_{-1}(\Sigma_0)),\nonumber
\end{align}
and the equality holds if and only if $\Sigma_0$ is some coordinate slice.
\end{theorem}
\begin{remark}
Notice that
\begin{align}
\mathcal{B}_{-1}(\Sigma)=\int_{\Sigma} u~ d\mu\nonumber
\end{align}
for any closed spacelike hypersurface in de Sitter space $\mathbb{S}^n_1$. Let $\Sigma_0$ be a closed spacelike convex hypersurface in de Sitter space $\mathbb{S}^n_1$, which meets further the pinching condition \eqref{pc}, then there exists a closed static convex hypersurface $\widetilde\Sigma_0$ in hyperbolic space $\mathbb{H}^{n+1}$. According to the dual relationship discovered in \cite{JS},
\begin{align}
\int_{\Sigma_0} u~ d\mu=\int_{\widetilde\Sigma_0} \widetilde{\lambda'}\widetilde{E}_n ~d\tilde{\mu}=\mathcal{B}_n(\widetilde{\Sigma}_0).\nonumber
\end{align}
On the other hand, due to the Minkowski formula, there holds
\begin{align}
\mathcal{B}_k(\Sigma_0)=\int_{\widetilde{\Sigma}_0}\widetilde{u}~\widetilde{E}_{n-k}~d\tilde{\mu}=\int_{\widetilde{\Sigma}_0}\widetilde{\lambda'}~\widetilde{E}_{n-k-1}~d\tilde{\mu}=\mathcal{B}_{n-k-1}(\widetilde{\Sigma}_0).\nonumber
\end{align}
In other words, once our assumption over the Heintze-Karcher inequality for closed spacelike mean convex hypersurface in de Sitter space can be removed, one particular sequence of geometric inequalities conjectured for closed static convex hypersurface in hyperbolic space $\mathbb{H}^{n+1}$ follows from Theorem \ref{MT2}.
\end{remark}
The rest of this paper is organized as follows. In Section \ref{Pre}, we recall fundamental formulae concerning the geometry of spacelike hypersurface in de Sitter space $\mathbb{S}^n_1$, and derive the uniform $C^0$ estimate along flow \eqref{f1}. In Section \ref{GE}, we derive for important geometric quantities the evolution equations along flow \eqref{f1} and obtain the uniform gradient estimate. In Section \ref{PPC}, we show the preservation of pinching condition \eqref{pc} along flow \eqref{f1}, applying the maximum principle for tensors. In Section \ref{PKC}, we show the preservation of $k$-convexity along flow \eqref{f1}, which implies the uniform curvature estimate, comparing the pinching condition \eqref{pc}. In Section \ref{GI}, the convergence result of flow \eqref{f1} is claimed and a sequence of Alexandrov-Fenchel inequalities is established.

\section{Preliminaries and uniform $C^0$ estimate}\label{Pre}
Let $X, Y$ be tangential vector fields on the evolving hypersurface along flow \eqref{f1},
\begin{align}
g(X, Y)=\bar g(x_*(X), x_*(Y)),\nonumber
\end{align}
defines the induced metric on the evolving hypersurface, $x_*$ is the corresponding push-forward map. In the sequel, we shall not tell a tangential vector field along the evolving hypersurface from its push-forward. Define the scalar-valued second fundamental form $h$ via the decomposition
\begin{align}\label{gf}
\overline\nabla_X Y=\nabla_X Y+h(X, Y)\nu,
\end{align}
where $\nabla$ is the Levi-Civita connection determined by $g$ and $\nu$ is the future-directed timelike unit normal vector field. Notice that \eqref{gf} is hugely different from its counterpart in Riemannian setting. As a consequence,
\begin{align}
h(X, Y)=-\bar g(\overline\nabla_X Y, \nu)=\bar g(Y, \overline\nabla_X\nu)= g(Y, \mathcal{W}(X)).\nonumber
\end{align}
Combining the conventions above, there holds the Gauss equation
\begin{align}
\bar R(W, Z, X, Y)=R(W, Z, X, Y)+h(Y,Z)h(X, W)-h(X, Z) h(Y, W)\nonumber
\end{align}
for any tangential vector fields $X, Y, Z, W$ along the evolving hypersurface, where
\begin{align}
\bar R (W, Z, X, Y)=\bar g(\overline\nabla_X\overline\nabla_Y Z, W)-\bar g(\overline\nabla_Y\overline\nabla_X Z, W)-\bar g(\overline\nabla_{[X,Y]}Z, W)\nonumber
\end{align}
is the curvature tensor of type $(0,4)$. Similarly, the Codazzi equation reads
\begin{align}
0=\bar R(\nu, Z, X, Y)=(\nabla_Y h)(Z, X)-(\nabla_X h)(Z,Y),\nonumber
\end{align}
where the covariant derivative
\begin{align}
(\nabla_X h)(Y,Z)=X(h(Y,Z))-h(\nabla_X Y , Z)-h(Y, \nabla_X Z).\nonumber
\end{align}
Define further the second order covariant derivative of $h$ as
\begin{align}
(\nabla^2 h)(Z, W; X, Y)=(\nabla_Y (\nabla_X h))(Z,W)-(\nabla_{(\nabla_Y X)} h)(Z, W).\nonumber
\end{align}
It is easy to check the following Ricci identity
\begin{align}
(\nabla^2 h)(Z, W; X, Y)-(\nabla^2 h)(Z, W; Y, X)=h(R(X, Y)Z, W)+h(Z, R(X, Y)W),\nonumber
\end{align}
where for any tangential vector field $U$ along the evolving hypersurface
\begin{align}
g (R(X, Y)Z, U)=R(U, Z, X, Y).\nonumber
\end{align}

For convenience, we express geometric quantities mentioned above under local coordinates and suppose the evolving hypersurface $\Sigma_t$ can be written locally as radial graphs
\begin{align}
\big\{(r(\theta, t),\theta)~|~ \theta\in\mathbb{S}^n\big\},\nonumber
\end{align}
where $r(\cdot, t): \mathbb{S}^n\rightarrow (0, +\infty)$ is a family of smooth functions. Now the induced structure $g$ on $\Sigma_t$ is locally determined by
\begin{align}
g_{ij}=\lambda^2\sigma_{ij}-r_i r_j.\nonumber
\end{align}
The evolving hypersurface $\Sigma_t$ is spacelike if and only if
\begin{align}
1-\lambda^{-2}|Dr|^2>0,\nonumber
\end{align}
where $D$ is the Levi-Civita connection determined by $\sigma$ and the norm is with respect to $\sigma$ as well. For convenience, let
\begin{align}
\upsilon=\sqrt{1-\lambda^{-2}|Dr|^2}\nonumber
\end{align}
whenever the evolving hypersurface is spacelike. Now the future-directed timelike unit normal vector field is
\begin{align}
\nu=\frac{1}{\upsilon}\big(\partial_r+\lambda^{-2}Dr\big).\nonumber
\end{align}
After a lengthy computation, it turns out that the second fundamental form can be expressed locally as
\begin{align}
h_{ij}=\frac{1}{\upsilon}\bigg(r_{, ij}+\lambda\lambda'\sigma_{ij}-2\frac{\lambda'}{\lambda}r_i r_j\bigg),\nonumber
\end{align}
where $D^2r=[r_{, ij}]$ is the Hessian with respect to $\sigma$. Or equivalently,
\begin{align}
h_{ij}=\upsilon\big(r_{; ij}+\lambda\lambda'\sigma_{ij}\big),\nonumber
\end{align}
with $\nabla^2r=[r_{; ij}]$ the Hessian with respect to the induced metric $g$. In particular,
\begin{align}\label{hor}
\Lambda_{; pq}=u h_{pq}-\lambda' g_{pq}.
\end{align}
As a consequence, noticing the divergence-free property of the tensor $E_l^{ij}$, there holds the Minkowski formula
\begin{align}
\int_{\Sigma} u E_l~ d\mu=\int_{\Sigma} \lambda' E_{l-1}~d\mu,
\end{align}
for every $l\in\{1,\cdots, n\}$ and every closed spacelike hypersurface $\Sigma$ in $\mathbb{S}^n_1$.

Let
\begin{align}
F=\frac{E_k}{E_{k-1}},\nonumber
\end{align}
then $F$ is obviously homogeneous of degree one. Besides, it is well known that $F$ is monotone increasing and concave in $\Gamma_k^+$, hence there hold the inequalities
\begin{align}
\frac{\partial F}{\partial \kappa_i}(\kappa) >0 ~\;~\text{and} ~\;~ \frac{\partial^2 F}{\partial\kappa_i\partial\kappa_j}\leq 0\nonumber
\end{align}
for all $i, j\in\{1,\cdots, n\}$ and every $\kappa=(\kappa_1,\cdots,\kappa_n)\in \Gamma_k^+$. Furthermore, from the Newton inequalities
\begin{align}
\frac{E_{k+1}}{E_k}\leq \frac{E_k}{E_{k-1}}\leq \frac{E_{k-1}}{E_{k-2}}
\end{align}
it yields the important inequalities
\begin{align}\label{BIFCF}
\sum_{i=1}^n\frac{\partial F}{\partial \kappa_i}(\kappa)\geq 1~\;~\text{and}~\;~ \sum_{i=1}^n \kappa_i^2\frac{\partial F}{\partial \kappa_i}(\kappa)\geq F^2
\end{align}
for every $\kappa=(\kappa_1,\cdots,\kappa_n)\in \Gamma_k^+$. Sometimes, it is more convenient to consider $F$ as smooth functions in entries of the induced structure and the second fundamental form on spacelike hypersurface. For simplicity, let
\begin{align}
F^{pq}=\frac{\partial F}{\partial h_{pq}} ~\;~ \text{and}~\;~  F^{pq,rs}=\frac{\partial^2 F}{\partial h_{pq}\partial h_{rs}}.
\end{align}

Using De Turck's trick, flow \eqref{f1} is equivalent to
\begin{align}\label{f2}
\frac{\partial}{\partial t} r=\frac{\lambda}{\lambda'}-\frac{\upsilon}{F}.
\end{align}
At the spatial maximum point of function $r$, $Dr$ vanishes, while $D^2r$ is non-positive definite. Hence the spatial maximum of function $r$, as a Lipschitz function in $t$, is monotone decreasing, where the homogeneity and monotonicity of $F$ are involved. Applying the same argument to the spatial minimum of function $r$, we obtain finally the uniform $C^0$ estimates along flow \eqref{f2}. More precisely, there exists some constant $C>1$ depending only on the initial value such that along flow \eqref{f2}
\begin{align}
C^{-1}\leq r \leq C.\nonumber
\end{align}
By definition, the support function
\begin{align}
u=\frac{\lambda}{\upsilon}\nonumber
\end{align}
is obviously no less than $\lambda$. As a consequence, the support function is uniformly bounded from below naturally along flow \eqref{f1}.
\section{Evolution equations and uniform gradient estimate}\label{GE}
From \eqref{ck}, it yields that
\begin{align}\label{gos}
\nabla u=\mathcal{W}(\nabla \Lambda).
\end{align}
Furthermore, the Hessian of the support function can be expressed locally as
\begin{align}\label{hos}
u_{;pq}=-\lambda' h_{pq}+g(\nabla\lambda', \nabla h_{pq})+u h_p^m h_{mq}.
\end{align}
Due to \eqref{hor}, \eqref{gos} and \eqref{hos}, there hold
\begin{align}\label{goT}
\nabla\Theta=(\mathcal{W}-\Theta\mathcal{I})(\nabla\ln\lambda'),
\end{align}
and
\begin{align}\label{hoT}
\nabla^2\Theta=&g(\nabla\ln\lambda', \nabla h )+\Theta h^2-(1+\Theta^2) h+\Theta g-\nabla(\ln \lambda')\otimes\nabla\Theta,
\end{align}
where
\begin{align}
[\nabla(\ln \lambda')\otimes\nabla\Theta] (X, Y)=\nabla_X(\ln \lambda')\otimes\nabla_Y\Theta+\nabla_Y(\ln \lambda')\otimes\nabla_X\Theta.\nonumber
\end{align}

For convenience, define the second order parabolic operator
\begin{align}
\mathscr{L}\cdot=\frac{\partial}{\partial t}\cdot-\frac{1}{F^2} F^{pq}\nabla^2_{pq}\cdot-g(\nabla\ln\lambda' ,  \nabla \cdot )\nonumber
\end{align}
along flow \eqref{f1}.

\begin{lemma}
Along flow \eqref{f1}, there hold the evolution equations
\begin{align}\label{eow}
\mathscr{L}\lambda'=\frac{\lambda^2}{\lambda'} -2\lambda'\frac{\Theta}{F}+\frac{\lambda'}{F^2}F^{pq} g_{pq},
\end{align}
\begin{align}\label{eos}
\mathscr{L} u=u\bigg(1-\frac{1}{F^2}F^{pq}h_p^m h_{mq}\bigg)-u||\nabla\ln\lambda'||^2,
\end{align}
and
\begin{align}\label{popc1}
\mathscr{L}\Theta&=\Theta\bigg(1-\frac{1}{F^2}F^{pq}h_p^m h_{mq}\bigg)-\Theta\bigg(\Theta - \frac{1}{F}\bigg)^2+\frac{\Theta}{F^2}(1-F^{pq}g_{pq})\nonumber\\
&+\frac{2\Theta}{F^2}F^{pq}(\ln\lambda')_{;p}(\ln\Theta)_{; q},
\end{align}
where the norm $||\cdot||$ is with respect to the induced metric on evolving hypersurface.
\end{lemma}
\begin{proof}
On the one hand, there holds
\begin{align}
\frac{\partial }{\partial t} \lambda'=\bar g(\overline\nabla\lambda', \mathcal{F} \nu) =u\mathcal{F},\nonumber
\end{align}
as well as
\begin{align}
\frac{\partial }{\partial t} u =\lambda'\mathcal{F}+g(\nabla\lambda',\nabla\mathcal{F}).\nonumber
\end{align}
Here the conformal Killing property \eqref{ck} and
\begin{align}
\frac{\partial}{\partial t}\nu=\nabla\mathcal{F}\nonumber
\end{align}
are used.
On the other hand, contracting \eqref{hor} and \eqref{hos} with $F^{pq}$ gives
\begin{align}
-\frac{1}{F^2}F^{pq} (\lambda')_{;pq}=-\frac{u}{F}+\frac{\lambda'}{F^2} F^{pq}g_{pq}\nonumber
\end{align}
and
\begin{align}
-\frac{1}{F^2}F^{pq} u_{;pq}=\frac{\lambda'}{F}-\frac{\lambda'}{F^2} g(\nabla\ln\lambda', \nabla F)-\Theta\frac{\lambda'}{F^2} F^{pq} h_p^m h_{mq}\nonumber
\end{align}
respectively. Besides, the identity
\begin{align}\label{ii}
||\nabla\ln\lambda'||^2=\Theta^2-\frac{\lambda^2}{\lambda'^2}
\end{align}
is also involved. According to the products rule for taking derivatives, the evolution equation of the function $\Theta$ follows directly.
\end{proof}
Now the uniform upper bound of the support function follows from the classical maximum principle. More precisely, the spatial maximum of the support function is monotone decreasing in $t$. In other words, there exists some constant $C>1$ depending only on the initial value such that along flow \eqref{f1}
\begin{align}
C^{-1}\leq u\leq C.\nonumber
\end{align}
And this implies the preservation of being spacelike along flow \eqref{f1}. In fact, the evolving hypersurface fails to be spacelike if and only if its support function blows up.

\section{preservation of pinching condition \eqref{pc} along flow \eqref{f1}}\label{PPC}

Let
\begin{align}
S^i_j=\Theta \delta^i_j- h^i_j,\nonumber
\end{align}
where $\delta^i_j$ is the Christoffel symbol. Notice that the pinching condition \eqref{pc} is equivalent to $S^i_j\geq 0$. In this section, we shall derive the evolution equation of the tensor $S^i_j$ first, and then apply the maximum principle for tensors developed in \cite{BA} to claim that the pinching condition is preserved along flow \eqref{f1}. More precisely, Andrews proved the following

\begin{proposition}
Let $S_{ij}$ be a smooth time-varying symmetric tensor field on a compact manifold $M$ without boundary, which satisfies
\begin{align}
\frac{\partial}{\partial t}S_{ij}=a^{kl}\nabla_k\nabla_lS_{ij}+u^k\nabla_kS_{ij}+N_{ij},\nonumber
\end{align}
where $a^{kl}$ and $u$ are smooth, $\nabla$ is a smooth symmetric connection which is possibly time-dependent, and $a^{kl}$ is positive definite everywhere. Suppose that
\begin{align}\label{MPFT}
N_{ij}v^iv^j+\sup_{\Lambda} 2a^{kl}(2\Lambda_k^p\nabla_l S_{ip}v^i-\Lambda_k^p\Lambda_l^q S_{pq})\geq 0
\end{align}
whenever $S_{ij}\geq 0$ and $S_{ij} v^j=0$. Here the supremum is taken over all $n\times n$ matrices. If $S_{ij}$ is positive definite everywhere on $M$ at time $t=0$, then it is positive on $M \times[0,T]$.
\end{proposition}

\begin{lemma}
Along flow \eqref{f1}, the Weingarten transformation evolves as
\begin{align}\label{popc2}
\mathscr{L} h^i_j&=\frac{2}{F}h^i_l h^l_j-\bigg(1+\Theta^2 +\frac{1}{F^2}F^{pq}h_p^mh_{mq}+\frac{1}{F^2}F^{pq}g_{pq}\bigg)h^i_j+2\Theta \delta^i_j\nonumber\\
&-\big[(\ln \lambda')_{; }^{\; i}  \Theta_{; j}+\Theta_{; }^{\; i}(\ln \lambda')_{; j}\big]+\frac{1}{F^2}F^{pq,rs}h_{pq;}^{\; \;  \; \; i}h_{rs; j}.
\end{align}
\end{lemma}
\begin{proof}
Firstly,
\begin{align}
\frac{\partial}{\partial t} h^i_j=\frac{\partial}{\partial t} (g^{il} h_{lj})=\mathcal{F}\delta^i_j+\mathcal{F}^{\; i}_{; \; \, j}-\mathcal{F} h^i_l h^l_j,
\end{align}
since
\begin{align}
\frac{\partial}{\partial t} g_{ij}=2\mathcal{F} h_{ij},\nonumber
\end{align}
while
\begin{align}
\frac{\partial}{\partial t} h_{ij}=\mathcal{F} g_{ij}+\mathcal{F}_{; ij}+\mathcal{F} h_i^l h_{lj}.\nonumber
\end{align}
Here the fact that $\mathbb{S}^n_1$ is of constant sectional curvature $1$ is invoked.

Secondly,
\begin{align}
\bigg(-\frac{1}{F}\bigg)^{\; i}_{; \; j}=\frac{1}{F^2}F^{pq} h^{\; ~\;~  ~\;~ i}_{pq; \; j}-\frac{2}{F^3}F^{pq, rs} h^{ ~\;~ ~\;~ i}_{pq; } h_{rs; j},\nonumber
\end{align}
while the Gauss equation, the Codazzi equation and the Ricci identity imply that
\begin{align}
F^{pq}h^{\; \; \; \; \, i}_{pq; \; j}=F^{pq}h^i_{j; pq}+ F \delta^i_j-F^{pq}h_p^mh_{mq}h^i_j- F^{pq}g_{pq} h^i_j+ F h_l^i h^l_j.\nonumber
\end{align}
Noticing that the Hessian of function $\Theta$ is given in \eqref{hoT}.
\end{proof}
Now it follows from \eqref{popc1} and \eqref{popc2} that the tensor $S^i_j$ evolves according to
\begin{align}
\mathscr{L} S^i_j&=-\frac{2}{F}S^i_l S^l_j-\bigg(1+\Theta^2-\frac{4\Theta}{F}+\frac{1}{F^2}F^{pq}h_p^mh_{mq}+\frac{1}{F^2}F^{pq}g_{pq}\bigg)S^i_j\nonumber\\
&+\frac{2}{F^2}F^{pq}(\ln\lambda')_{;p}\Theta_{; q}\delta^i_j+\big[(\ln \lambda')_{; }^{\; i}  \Theta_{; j}+\Theta_{; }^{\; i}(\ln \lambda')_{; j}\big]-\frac{1}{F^2}F^{pq,rs}h_{pq;}^{\; \;  \; \; i}h_{rs; j}.\nonumber
\end{align}

Let $(t^*, p^*)$ be the point where $S^i_j$ has an unit null vector $v$. By continuity, we may assume that the principal curvatures are mutually distinct and in decreasing order at $(t^*, p^*)$. Now the null vector condition $S^i_j v^j=0$ implies that $v=e_1$ and $S^1_1=\Theta-\kappa_1=0$ at $(t^*, p^*)$. It remains to show
\begin{align}\label{MPFT-1}
Q:=&\frac{2}{F^2}\sum_{l=2}^nF^{ll}(\ln\lambda')_{;l}\Theta_{; l}-\frac{1}{F^2}F^{pq,rs}h_{pq;1}h_{rs; 1}\nonumber\\
&+\sup_{\Gamma}\frac{2}{F^2} F^{kl}(2\Gamma_k^p\nabla_lS_{1p}-\Gamma_k^p\Gamma_l^q S_{pq})\geq 0.
\end{align}
By assumption, $S_{11}=0$ and $\nabla S_{11}=0$ at the considered point. We have
\begin{align}
F^{kl}(2\Gamma_k^p\nabla_lS_{1p}-\Gamma_k^p\Gamma_l^q S_{pq})&=\sum_{l=1}^n\sum_{p=2}^nF^{ll}[2\Gamma_l^p\nabla_lS_{1p}-(\Gamma_l^p)^2 S_{pp}]\nonumber\\
&=\sum_{l=1}^n\sum_{p=2}^nF^{ll}\bigg[\frac{(\nabla_l S_{1p})^2}{S_{pp}}-\bigg(\Gamma_l^p-\frac{\nabla_l S_{1p}}{S_{pp}}\bigg)^2 S_{pp}\bigg].\nonumber
\end{align}
Hence the supremum in \eqref{MPFT-1} is achieved by choosing $\Gamma_l^1=0$ for all $l\in\{1,\cdots,n\}$ and
\begin{align}
\Gamma_l^p=\frac{\nabla_l S_{1p}}{S_{pp}}~\;~ \text{for}~\;~ l,p\in\{2,\cdots,n\}.\nonumber
\end{align}
It follows that
\begin{align}
Q&=\frac{2}{F^2}\sum_{l=2}^nF^{ll}(\ln\lambda')_{;l}\Theta_{; l}-\frac{1}{F^2}F^{pq,rs}h_{pq;1}h_{rs; 1}+\frac{2}{F^2}\sum_{l=1}^n\sum_{p=2}^nF^{ll}\frac{(\nabla_l S_{1p})^2}{S_{pp}}\nonumber\\
&\geq \frac{2}{F^2}\sum_{l=2}^nF^{ll}(\ln\lambda')_{;l}\Theta_{; l}-\frac{2}{F^2}\sum_{p=2}^n\frac{F^{pp}-F^{11}}{\kappa_p-\kappa_1}(\nabla_1h_{1p})^2+\frac{2}{F^2}\sum_{p=2}^nF^{11}\frac{(\nabla_1 h_{1p})^2}{S_{pp}}\nonumber\\
&=\frac{2}{F^2}\sum_{l=2}^nF^{ll}(\ln\lambda')_{;l}\Theta_{; l}+\frac{2}{F^2}\sum_{p=2}^n\frac{F^{pp}}{S_{pp}}(\nabla_p h_{11})^2\nonumber\\
&=\frac{2}{F^2}\sum_{l=2}^nF^{ll}(\ln\lambda')_{;l}\Theta_{; l}-\frac{2}{F^2}\sum_{p=2}^nF^{pp}(\ln\lambda')_{;p}\Theta_{;p}=0,\nonumber
\end{align}
where the concavity of $F$, the Codazzi equation and the identity
\begin{align}
\Theta_{;p}=(\kappa_p-\Theta)(\ln\lambda')_{;p}=-S_{pp}(\ln\lambda')_{;p}\nonumber
\end{align}
are involved. In particular, there holds for any symmetric matrix $\eta$ the identity
\begin{align}
F^{pq,rs} \eta_{pq} \eta_{rs}=\sum_{i,j}\frac{\partial^2 F}{\partial\kappa_i\partial\kappa_j} \eta_{ii}\eta_{jj}+\sum_{i,j}\frac{\frac{\partial F}{\partial\kappa_i}-\frac{\partial F}{\partial\kappa_i}}{\kappa_i-\kappa_j}\eta_{ij}^2.\nonumber
\end{align}
\section{preservation of $k$-convexity}\label{PKC}
\begin{lemma}
Along flow \eqref{f1}, the curvature function $F$ evolves according to
\begin{align}\label{eeof}
\mathscr{L} F&=2\Theta F^{pq} g_{pq}-\frac{1}{F} F^{pq} g_{pq}-(1+\Theta^2) F+\frac{1}{F}F^{pq} h_p^m h_{mq},\nonumber\\
&-2F^{pq} (\ln\lambda')_{;p} \Theta_{;q}-2\frac{1}{F^3}F^{pq} F_{;p} F_{;q}.
\end{align}
\end{lemma}
\begin{proof}
Substituting \eqref{hoT} into
\begin{align}
\frac{\partial }{\partial t} F=\mathcal{F}F^{pq} g_{pq}+F^{pq}\Theta_{;pq}+\frac{1}{F^2} F^{pq}F_{;pq}-2\frac{1}{F^3}F^{pq}F_{;p} F_{;q}-\mathcal{F}F^{pq}h_p^m h_{mq},\nonumber
\end{align}
then the desired evolution equation follows.
\end{proof}

We consider the auxiliary function
\begin{align}
\omega=-\ln (\lambda'-\delta)-\ln F,\nonumber
\end{align}
where
\begin{align}
\delta=\frac{1}{2}\min \lambda'.\nonumber
\end{align}
From \eqref{eow} and \eqref{eeof}, it yields directly the evolution equation of $\omega$. More precisely,
\begin{align}
\mathscr{L} \omega&=\frac{\lambda'}{\lambda'-\delta}||\nabla\ln\lambda'||^2-\frac{\delta}{\lambda'-\delta}\bigg(\frac{1}{F}-\Theta\bigg)^2+\frac{2}{F}F^{pq} (\ln\lambda')_{;p} \Theta_{;q}\nonumber\\
&+\frac{\delta}{\lambda'-\delta}\frac{1}{F^2}(1-F^{pq} g_{pq})+\frac{2\Theta}{F}(1-F^{pq} g_{pq})+1-\frac{1}{F^2}F^{pq} h_p^m h_{mq}\nonumber\\
&+\frac{1}{F^2}F^{pq} (\ln F)_{;p} (\ln F)_{;q}-\frac{1}{F^2}F^{pq} [\ln(\lambda'-\delta)]_{;p} [\ln(\lambda'-\delta)]_{;q},\nonumber
\end{align}
where \eqref{ii} is used again. Hence there holds the inequality
\begin{align}
0& \leq -\frac{\delta}{\lambda'-\delta}\bigg(\frac{1}{F}-\Theta\bigg)^2+\frac{\lambda'}{\lambda'-\delta}||\nabla\ln \lambda'||^2\nonumber
\end{align}
at the point where $\omega$ attains its maximum, where the critical point condition, the pinching condition and inequalities \eqref{BIFCF} are used. Hence there exists some constant $C>0$ depending only on the initial hypersurface such that
\begin{align}
F\geq C\nonumber
\end{align}
along flow \eqref{f1}. As a consequence, the evolving hypersurface along flow \eqref{f1} preserves to be $k$-convex if the initial one is.

From the Newton inequality and the pinching condition \eqref{pc}, it follows the uniform curvature estimate along flow \eqref{f1}. According to the standard regularity theory, flow \eqref{f1} exists for all positive time.

\section{Convergence and Alexandrov-Fenchel type inequalities}\label{GI}
It is easy to check that
\begin{align}\label{vf}
\frac{d}{dt} \mathcal{B}_k(\Sigma_t)&=(1+k)\int_{\Sigma_t} u E_k \mathcal{F}\, d\mu+(n-k)\int_{\Sigma_t}\lambda' E_{k+1} \mathcal{F}\, d\mu\nonumber\\
&\geq \frac{1+k}{k}\int_{\Sigma_t} \Theta E_k^{ij}(\lambda')_{;ij}\,d\mu\nonumber\\
&=\frac{1+k}{k}\int_{\Sigma_t} [(\Theta\mathcal{I}-\mathcal{W}) (\nabla\lambda')]_{;j} E_k^{ij}(\ln \lambda')_{;i}\,d\mu\nonumber\\
&\geq 0
\end{align}
along flow \eqref{f1}, where the Newton inequality, the divergence theorem and the pinching condition \eqref{pc} are used. Notice that one necessary condition such that the equality in \eqref{vf} holds is that $\Sigma_t$ is umbilical. On the other hand, due to the a priori estimates, $\mathcal{B}_k(\Sigma_t)$ is uniformly bounded along flow \eqref{f1}, which forces the limiting hypersurface to be umbilical, or equivalently to be some coordinate slice.

Due to the Heintze-Karcher inequality and the Newton inequality, the weighted curvature integral $\mathcal{B}_{-1}(\Sigma_t)$ satisfies
\begin{align}
\frac{d}{dt}\mathcal{B}_{-1}(\Sigma_t)=(n+1)\int_{\Sigma_t}\bigg( u-\frac{\lambda'}{F} \bigg)d\mu(n+1)\int_{\Sigma_t}\bigg( u-\frac{\lambda'}{E_1} \bigg)d\mu\leq 0\nonumber
\end{align}
along flow \eqref{f1}.

Let $\Sigma_\infty=\{r_\infty\}\times \mathbb{S}^n$ be the limiting hypersurface along flow \eqref{f1}, then
\begin{align}
\mathcal{B}_{-1}(\Sigma_\infty)=\omega_n \lambda^{n+1}(r_\infty):=\varphi_{-1}(r_\infty),\nonumber
\end{align}
while
\begin{align}
\mathcal{B}_{k}(\Sigma_\infty)=\omega_n \lambda^{n-k}(r_\infty)(\lambda')^{k+1}(r_\infty):=\varphi_{k}(r_\infty).\nonumber
\end{align}
According to their monotonicity,
\begin{align}
\mathcal{B}_k(\Sigma_0)\leq \mathcal{B}_k(\Sigma_\infty)=\varphi_{k}(r_\infty)=\varphi_{k}\circ\varphi_{-1}^{-1}(\varphi_{-1}(r_\infty))\leq \varphi_{k}\circ\varphi_{-1}^{-1}(\mathcal{B}_{-1}(\Sigma_0)),\nonumber
\end{align}
since both $\varphi_{-1}$ and $\varphi_k$ are strictly monotone increasing.
\section{Acknowledgement}
This work has been supported by Natural Science Foundation of Shandong Province (Grant No. ZR202111230079).

\end{document}